\documentclass[11pt]{amsart}
\usepackage{amsmath,amscd} 
\usepackage{latexsym,amsmath,amssymb,mathrsfs,setspace,enumerate}
\usepackage[T1]{fontenc}
\usepackage[utf8]{inputenc}
\usepackage{lmodern} 
\usepackage{amssymb}
\usepackage{tikz}
\usepackage{pst-node}
\usepackage{tikz-cd}
\newtheorem{lemma}{Lemma} 
 
\newtheorem{thm}{Theorem}

\newtheorem{exam}{Example}
\newtheorem{cor}{Corollary} 
\newtheorem{dfn}{Definition} 
\newcommand{\bl}{\begin{lemma}}
\newcommand{\el}{\end{lemma}}
\newcommand{\bt}{\begin{thm}}
\newcommand{\et}{\end{thm}}
\newcommand{\bc}{\begin{cor}}
\newcommand{\ec}{\end{cor}}
\newcommand{\bd}{\begin{defn}}
\newcommand{\ed}{\end{defn}}
\newcommand{\bp}{\begin{proof}}
\newcommand{\ep}{\end{proof}}
\newcommand{\be}{\begin{exam}}
\newcommand{\ee}{\end{exam}}
\newcommand{\oml}{\omega^l}
\newcommand{\omr}{\omega^r}
\def \ni{\noindent}
\def \ml{\mathrel}

\title{Biordered Sets of Lattices and Homogeneous Basis}
\author { P. G. Romeo$^1$ and Akhila. R$^2$ }
\address{$^1$ Professor, Department of Mathematics, Cochin University of Science and Technology, Kochi, Kerala, INDIA.}
\address{$^2$ Assistant Professor, Department of Mathematics, Amrita School of Arts and Sciences, Kochi, Amrita Vishwa Vidyapeetham, INDIA}

\email{$^1 romeo_-parackal@yahoo.com\,,  ^2 akhila.ravikumar@gmail.com, $}
\date{}

\begin{document}

\begin{abstract} 

In this paper, we discuss the properties of the biordered set $E_{P(L)}$ obtained from the complemented modular lattice $L$, defines an operation $\oplus$ using the sandwich elements of 
the biordered set. Further, we describe the biordered subset $E_{P(L)}^0$ of $E_{P(L)}$ satisfying certain conditions, so that the complemented modular lattice admits a homogeneous basis. Finally, analogous to von Neumann's coordinatization theorem we describe the coordinatization theorem for complemented modular lattice using the biordered set of idempotents $E_{P(L)}^0$.

\end{abstract}
\maketitle
\section{Introduction}

The concept of regular rings was  introduced by John von Neumann in a paper "On Regular Rings" in 1936. He used regular rings as an algebraic tool for studying certain lattice of projections on algebra of operators on a Hilbert space. A lattice $L$ is said to be coordinatized by a regular ring $R$, if $L$ is isomorphic to the principal left[right] ideals of the regular ring.
von Neumann proved that every complemented modular lattice with order greater than or equal to $4$ is coordinalizable (see \cite{von Neu}).  
\ni  

 The multiplicative reduct of a regular ring is a regular semigroup and thus it is obvious that the study of regular semigroups play a singificant role in the study of regular rings.  
 In order to study the structure of regular semigroup, in his celebrated 1973, paper K.S.S. Nambooripad introduced the concept of a biordered sets to study the structure of idempotents of a regular semigroup (cf.\cite{Namb}). He defined a biordered set as an order structure to represent the set of idempotents of a semigroup and identified a partial binary operation on 
 the set of all idempotents $E(S)$ of the semigroup $S$ arising from the underlined binary operation in $S$ and defined two quasiorders $\omega^r$ and $\omega^l$
on the set $E(S)$ satisfying certain axioms (see the definition below) which is abstracted as a biordered set.\\

In \cite{Akhi2}  we extend the biordered set approch from regular semigroups to regular rings by explicitly desceribing the structure of the multiplicative idempotents $E_R$ of a 
regular ring $R$ with two quasi orders $\omega^l$ and $\omega^r$ as a bounded and complemented biorderd set . In \cite{Akhi}, we describe the biorder ideals generated by 
elements in $E_R$ as 
$$\{\omega^l(e): \,e\in E_R\}\,\,\text{and}\,\,\{\omega^r(e): \, e\in E_R\}$$
\ni
and are the complemented modular lattice $\Omega_L$ and $\Omega_R$ respectively. Further we also discuss certain interesting propertes of the complemented modular lattice such as perspectivity and independence of the lattce and obtained the necessary and suffient condition at which the complementd modular lattice is of order $n$.
\ni 

The converse problem of describing the structure which coordinatise a given complemented modular lattice is attempted by Pastjin in the case of a strongly regular baer 
semigroup (see cf.\cite{pastjin}). He defined the normal mappings on a complemented modular lattice $L$ using complementary pairs and it is shown that these normal mappings form a semigroup 
$P(L)$ and the set of idempotents $E_{P(L)}$ of $P(L)$ is the biordered set of the complemented modular lattice and $L$ can be coordinatized by $P(L)$. 
\ni 

Here we extend Pastjin's approach to obtain a regular ring coordinatizing the complemented modular lattice $L$. 
For, we note that the biordered set $E_{P(L)}$ obtained from the complemented modular lattice $L$ is a biordered set with least and greatest element, each element having an inverse and 
for elements satisfying $v_i\leq n_j$ in $L$, an operation in $E_{P(L)}$ can be defined using the sandwich sets and we observe some properties of this operation such as cancellation, inverse. 
It is also seen that on some conditions in $E_{P(L)}$, the complemented modular lattice $L$ admits a homogeneous basis and this lattice $L$ corresponds to the lattice $\Omega_L$ the lattice of all $\omega^l$-ideals of a ring $R$ (cf.\cite{Akhi}).

\section{Preliminaries}
Let $S$ be a regular semigroup, and $E(S)$ the set of idempotents of $S$. The relation $\omega^l[\omega^r]$ on $E(S)$ is defined by $e\mathrel \omega^l f[e\mathrel \omr f]$ if and only 
if $ef = e[fe = e]$ is a quasiorder on $E(S)$; the sets,
$$\omega^r(e) = \left\{f: ef = f\right\}, ~\ \omega^l(e) = \left\{f:fe = f\right\}$$ 
are principal ideals and the relations $\mathcal R=\omega^r\cap(\omega^r)^{-1},\, \mathcal L=\omega^l\cap(\omega^l)^{-1}$, and $\omega(e) = \omega^r(e)\cap \omega^l(e)$
for all $e\in S$ are equivalences and partial order respectively.

\begin{dfn}
Let $E$ be a partial algebra. Then $E$ is a biordered set if the following axioms and their duals hold:
\begin{enumerate} 
\item  $\omega^r$ and $\omega^l$ are quasi orders on $E$ and $$D_E=(\omega^r\cup\omega^l)\cup(\omega^r\cup\omega^l)^{-1}$$ 
\item  $f\in \omega^r(e)\Rightarrow f\mathcal R fe\omega e$ \item $g\omega^l f\,\,\mbox{and}\,\, f,g\in \omega^r(e)\Rightarrow ge\omega^l fe.$ 
\item  $g\omega^r f\omega^r e\Rightarrow gf=(ge)f$ 
\item  $g\omega^l f\, \mbox{and} \,f,g\in \omega^r(e)\Rightarrow (fg)e=(fe)(ge).$ \\
Let $\mathcal M(e,f)$ denote the quasi ordered set $(\omega^l(e)\cap\omega^r(f),<)$ where $<$ is defined by 
$g< h\Leftrightarrow eg\mathrel{\omega^r} eh,\, \mbox{and}\,gf\mathrel{\omega^l} hf.$ Then the set 
$$ S(e,f)=\{h\in M(e,f): g< h\,\,\mbox{for all} \,\, g\in M(e,f)\}$$ is called the sandwich set of $e$ and $f$.
\item $f,g \in \omega^r(e)\Rightarrow S(f,g)e=S(fe,ge)$ 
\end{enumerate}
\end{dfn}

We shall often write
$E =\langle E,\omega^l,\omega^r\rangle$ to mean that $E$ is a biordered set with quasi-orders $\omega^l, \,\,\omega^r.$
 The relation $\omega$ defined is a partial order and 
$$\omega\cap (\omega)^{-1}\subset \omega^r\cap (\omega^l)^{-1}= 1_E.$$
The biordered set $E$ is said to be regular if
$S(e,f)\neq \emptyset \,\,\mbox{for all}\,\, e,f\in E$. 
Regular semigroups which determine the same biordered set will be called biorder isomorphic.

\begin{dfn}
Let $e$ and $f$ be idempotents in a semigroup $S$. By an $E$-sequence
from $e$ to $f$, we mean a finite sequence $e_0 = e, e_1, e_2, . . . , e_{n-1}, e_n = f$ of
idempotents such that $e_{i-1}(\mathcal{L}\cup\mathcal{R})e_i$ for $i = 1, 2, . . . , n$; $n$
is called the length of the $E$-sequence. 
\end{dfn}
If there exists an $E$-sequence from $e$
to $f$, $d(e, f)$ is the length of the shortest $E$-sequence from $e$
to $f$; and $d(e, e) = 1$. If there is no $E$-sequence from $e$ to $f$, we
define $d(e, f) = 0$.

\begin{dfn}(cf.\cite{Grill1})
If $P$ is a partially ordered set and $\Phi: P\rightarrow P$ is an isotone(order preserving) mapping, the $\Phi$ will be called normal if
	\begin{enumerate}
		\item $im\Phi$ is a principal ideal of $P$ and
		\item whenever $x\Phi = y$, then there exists some $z\leq x$ such that $\Phi$ maps the principal ideal $P(z)$ isomorphically onto the principal ideal $P(y)$.
	\end{enumerate}
	\end{dfn}
	
\begin{dfn}(cf. \cite{Grill1})
The partially ordered set $P$ will be called regular if for every $e \in P$, $P(e) = im\Phi$ for some normal mapping $\Phi: P\rightarrow P$ with $\Phi^2 = \Phi$ .
\end{dfn}
If $P$ is a partially ordered set, then it is easy to see that the set $S(P)[S^*(P)]$ of normal mappings of $P$ into itself, considered as left [right] operators form a regular semigroup.
\begin{dfn}
A lattice is a partially ordered set in which each pair of elements has the least upper bound and the greatest lower bound. If $a$ and $b$ are elements of a lattice, we denote their greatest lower bound
(meet) and least upper bound (join) by $a\wedge b$ and $a\vee b$ respectively. It is
easy to see that $a\vee b$ and $a\wedge b$ are unique. 
\end{dfn}

\begin{dfn}
A lattice is called modular (or a Dedekind lattice) if 
$$(a \vee b) \wedge c = a \vee (b \wedge c)\,\, \text{for all}\,\, a \leq c.$$ 
\end{dfn}
A lattice is bounded if it has both a maximum element and a minimum element, we use the symbols 0 and 1 to denote the minimum element and maximum element of a lattice. A bounded lattice $L$ is said to be complemented if for each
element $a$ of $L$, there exists at least one element $b$ such that $a \vee b = 1$ and
$a \wedge b = 0$.
The element $b$ is referred to as a complement of $a$. It is quite possible for an element of a complemented lattice to have many complements. 
An element $x$ is called a complement of $a$ in $b$ if $a\vee x = b$ and $a \wedge x = 0.$

\begin{dfn}
The elements $x_1, x_2, \cdots x_n$ of a lattice are called independent if
$$(x_1 \vee \cdots \vee x_{i-1} \vee x_{i+1} \vee \cdots \vee x_n) \wedge x_i = 0$$
for every $i$.
\end{dfn}
\begin{dfn}
Two elements $a$ and $b$ in a lattice $L$ are said to be perspective (in symbols $a\sim b$) if there exists $x$ in $L$ such that
$a \vee x = b \vee x ,a \wedge x = b \wedge x = 0.$
Such an element $x$ is called an axis of perspective.
\end{dfn}

\begin{dfn}
Let $L$ be a complemented, modular lattice with zero $0$ and unit $1$. By a basis of $L$ is meant a system $(a_i\colon i = 1,2,\hdots n)$ of $n$ elements of $L$ such that
$$(a_i; i = 1,2, \hdots n) {\rm \,\, are \, independent,\,\,} \,\,\,a_1\vee a_2\vee \hdots \vee a_n = 1$$
A basis is homogeneous if its elements are pairwise perspective.
$$a_i\sim a_j (i,j = 1,2, \hdots n)$$
The number of $n$ elements in a basis is called the order of the basis.
\end{dfn}
\begin{dfn}
A complemented modular lattice $L$ is said to have order $m$ in case it has a homogeneous basis of order $m$.
\end{dfn}

\section{Biordered set of the complemented modular lattice and its Homogeneous basis }
Let $L$ be a complemented modular lattice and $(n; v)$ be any pair of complementary elements of $L$\\
Let
$(n;v):L\rightarrow L,$ be the map defined by
$$ x\rightarrow v\wedge(n\vee x) \mbox{ for all } x \mbox{ in } L$$ 
and $(n;v)':L\rightarrow L$ defined by 
$$x\rightarrow n\vee(v\wedge x) \mbox{ for all }  x \mbox { in } L$$
are idempotent order preserving normal mappings of $L$ onto $[0,v]$ of $(L, \wedge)$ acts on $L$ as right operators denoted by $P(L)$ the subsemigroup of $S^*(L)$ which is generated by these idempotent normal mappings $(n;v), n, v \in L$ and the mappings $(n;v)'$ which are order preserving idempotent normal mapping of $L$ onto the principal ideal 
$[1,n]$ of $(L,\vee)$ into itself. Letting $(n;v)'$ act on $L$ as left operators, denote by $P(L)'$ the subsemigroup of $S(L)$ which is generated by these idempotent normal mappings $(n; v)', n, v \in L$.\\

Let $$E_{P(L)} = \left\{(n;v): n,v \in L\,\, n\vee v = 1, n\wedge v = 0\right\}$$and
	$$E_{P(L)'} = \left\{(n;v)': n,v \in L\,\, n\vee v = 1, n\wedge v = 0\right\}$$
	we refer to the elements $(n;v) [(n;v)']$ as idempotent generators of $P(L) [P(L)]'$.

	\bt[cf.\cite{pastjin}, Theorem 1]
	\begin{enumerate}
	Let L be a complemented modular lattice. Then
		\item $P(L)$ is a regular subsemigroup of $S^*(L)$ and\\
		$E_{P(L)} = \left\{(n;v): n,v \in L\,\, n\vee v = 1, n\wedge v = 0\right\}.$
		\item In $(E_{P(L)}, \omega^l, \omega^r)$ we have
		$$(n_1; v_1)\ml\omega^l (n_2; v_2) \Longleftrightarrow v_1\leq v_2 \mbox{ in } L$$
		and then $$(n_2; v_2)(n_1;v_1) = (n_2\vee(v_2\wedge n_1); v_1);$$
		we have
		$$(n_1; v_1)\ml\omega^r (n_2; v_2) \Longleftrightarrow n_2\leq n_1 \mbox{ in } L$$
		and then
		$$(n_1; v_1)(n_2;v_2) = (n_1; v_2\wedge(n_2\vee v_1))$$
		\item Let $(n_1; v_1)$ and $(n_2; v_2)$ be any idempotent of $P(L)$. Let $n$ be any complement of $v_1\vee n_2$ in $[n_2, 1]$; let $v$ be any complement of $v_1\wedge n_2$ in $[0 ,v_1]$; the $n$ and $v$ are complementary in $L$ and $(n;v)$ is an element in the sandwich set $S((n_1; v_1),(n_2; v_2))$. Conversely, any element in the sandwich set $S((n_1; v_1),(n_2; v_2))$ can be obtained in this way.	
	\end{enumerate}
	\et
	The above theorem provides a biordered set $E_{P(L)}$ from a complemented modular lattice $L$. The zero of $P(L)$ is $(1;0)$ and the identity is $(0;1)$. Obviously $(1;0)$ and $(0;1)$ are in $E_{P(L)}$ and for any $(n; v)$ in biordered set $E_{P(L)}$, $(n; v) \ml\omega (0;1)$ and $(1; 0)\ml\omega(n; v)$. 
	
\bl\label{4a}
	Let $(n_1;v_1)$, $(n_2;v_2)\in E_{P(L)}$, then
  \begin{enumerate}
	 \item $S((n_1;v_1),(n_2;v_2)) = S((n_2;v_2),(n_1;v_1)) = \left\{(1;0)\right\}$ if and only if  $v_1\leq n_2$ and $v_2\leq n_1$.
	\item  For $v_1\leq n_2$ and $v_2\leq n_1$, $(v_1\vee v_2; n_2\wedge n_1)$ is the unique element in $S((v_1; n_1),(v_2; n_2))$ and $S((v_2; n_2),(v_1; n_1))$. 
	\end{enumerate}
	\el
	\bp
	(1) $(n_1;v_1), (n_2;v_2)\in E_{P(L)},\,v_1\leq n_2$ and $(n;v)\in S((n_1;	v_1),(n_2;v_2))$. Then by definition of sandwich set $n$ is a complement of $n_2$ in $[n_2, 1]$ and $v$ is a complement of $v_1$ in $[0,v_1]$. Since, $(n; v)\ml \omega^l(n_1; v_1)$ and $(n; v)\ml\omega^r (n_2; v_2)$, it follows that $n_2\leq n$ and $n\vee n_2 = 1$ implies $n = 1$ and $v\leq v_1$ and $v\wedge v_1 = 0$ implies $v = 0$. Thus  $S((n_1;	v_1),(n_2;v_2)) = \left\{(1;0)\right\}$. Similarly, $S((n_2; v_2),(n_1;v_1)) = \left\{(1;0)\right\}$ when $v_2\leq n_1$. The converse follows immediately.\\
	
(2) For $v_1\leq n_2$ and $v_2 \leq n_1$, by definition,
	$$(v_1;n_1)(v_2;n_2) = (v_1\vee(n_1\wedge v_2); n_2\wedge(v_2\vee n_1)) = (v_1\vee v_2; n_2\wedge n_1)$$
	and $$(v_2;n_2)(v_1;n_1) = (v_2\vee(n_2\wedge v_1);n_1\wedge(v_1\vee n_2)) = (v_2\vee v_1; n_1\wedge n_2)$$ 
	thus $(v_1;n_1)(v_2;n_2) = (v_2;n_2)(v_1;n_1)$ and is easy to see that $(v_1\vee v_2)$ is a complement of $v_2\vee n_1$ in $[v_2, 1]$ and $(n_2\wedge n_1)$ is a complement of 
	$v_2\wedge n_1$ in $[0,n_1]$. Thus $(v_1\vee v_2; n_2\wedge n_1) 	\in S((v_1; n_1),(v_2; n_2))$. Similarly, $(v_1\vee v_2)$ is a complement of $v_1\vee n_2$ in $[v_1, 1]$ and 
	$(n_2\wedge n_1)$ is a complement of $v_1\wedge n_2$ in $[0,n_2]$. Therefore, $(v_1\vee v_2; n_2\wedge n_1)\in S((v_2; n_2),(v_1; n_1))$.\\
	
	It remains to prove the uniqueness of this element, suppose there exists another element 
	$(a;a')\in S((v_1; n_1)(v_2; n_2)) \cap S((v_2; n_2)(v_1; n_1))$. Then $v_1\leq a, a'\leq n_1, v_2\leq a, a'\leq n_2$ and from the definition of sandwich set it can be seen that $a\vee n_1 = 1, a'\wedge v_2 = 0, a\vee n_2 = 1, a'\wedge v_1 = 0$ and $a\wedge n_1 \leq v_2, n_1 \leq a'\vee v_2, a\wedge n_2 \leq v_1, n_2\leq a'\vee v_1$. Thus $a = v_1\vee v_2 $ and $a' = n_1\wedge n_2$. and $(a;a') = (v_1;n_1)(v_2;n_2)$ is unique.
	\ep
	
	Thus for each $(n;v) \in E_{P(L)}$ there exists an element $(v;n) \in E_{P(L)}$ such that $(n;v)(v;n) = (v;n)(n;v) = (1;0)$ and call this element $(v;n)$ the inverse of $(n;v)$. Also 
		\begin{itemize}
		\item $(n_1;v_1)\mathrel{\omega^l}(n_2;v_2)\Longleftrightarrow (v_2;n_2)\mathrel{\omega^r}(v_1;n_1)$.
		\item $v_1\leq n_2 \Longleftrightarrow S((n_1;v_1)(n_2;v_2)) = (1;0)$
	  \end{itemize}
		
	from here onwards we consider the biordered subset of $E_{P(L)}$ satisfying $v_i\leq n_j$ for all $(n_i; v_i)$, $(n_j; v_j)$ and $i\neq j$. Clearly in this biordered subset 
	$(v_i; n_i)(v_j; n_j) = (v_j; n_j)(v_i; n_i) = (v_i\vee v_j; n_j\wedge n_i)$. Define $$(n_i; v_i)\oplus (n_j;v_j) = (n_i\wedge n_j; v_i\vee v_j)$$  .

\bl
Consider the biordered subset of $E_{P(L)}$ with $v_i \leq n_j$ for $i\neq j$ and let $(p;q)$ denotes $(n_i;v_i)\oplus(n_j;v_j)$. Then we have the following:
\begin{enumerate}
	\item $(n_i;v_i), (n_j;v_j)\ml \in\omega((p;q))$
	\item If $(r;s) \in E_{P(L)}$ with $(n_i;v_i), (n_j;v_j)\in\,\mathrel{\omega^l}((r;s))$ then $(p;q)\in\,\mathrel{\omega^l}((r;s))$.
	\item If $(r;s) \in E_{P(L)}$ with $(n_i;v_i), (n_j; v_j) \in\,\mathrel{\omega^r}((r;s))$, then $(p;q)\in\,\mathrel{\omega^r}((r;s))$.
\end{enumerate}
\el
\bp (1) Note that $(p;q) \in S((n_i;v_i)(n_j;v_j))\cap S((n_j;v_j),(n_i;v_i))$. Therefore, $(q;p)\mathrel{\omega^l}(v_i;n_i), (q;p)\mathrel{\omega^r}(v_j;n_j), (q;p)\mathrel{\omega^l}(v_j;n_j)$ and $(q;p)\mathrel{\omega^r}(v_i;n_i)$. Thus $p\leq n_i, p\leq n_j, v_j\leq q, v_i\leq q$ and $(n_i;v_i)\mathrel{\omega}(p;q)$ and $(n_j; v_j)\ml\omega (p;q)$.\\

(2) Let $(r;s) \in E_{P(L)}$ with $(n_i;v_i)\mathrel{\omega^l}(r;s)$ and $(n_j;v_j)\mathrel{\omega^l}(r;s)$, then $v_i\leq s$ and $v_j\leq s$. Then as seen above in Lemma \ref{4a}.(2),
$$(q;p) = (v_i;n_i)(v_j;n_j).$$
Thus,
$v_j\leq s$ implies $(s;r) \mathrel{\omega^r}(v_j;n_j)$, that is $(v_j;n_j)(s;r) = (s;r)$. Similarly,
$v_i\leq s$ implies $(s;r) \mathrel{\omega^r}(v_i;n_i)$, that is $(v_i;n_i)(s;r) = (s;r)$.\\
\noindent
Therefore, $$(q;p)(s;r) =(v_i;n_i)(v_j;n_j)(s;r) = (s;r)$$ that is $(s;r)\mathrel{\omega^r}(q;p)$ also $(p;q)\mathrel{\omega^l}(r;s)$. In a similar manner $(3)$ follows.
\ep

The next lemma shows that the opertion $\oplus$ is cancellative.
\bl
Let $(n_i; v_i),(n_j; v_j), (n_k; v_k) \in E_{P(L)}$ with $v_i \leq n_j, v_j\leq n_i$ and $v_i \leq n_k , v_k \leq n_i$ for $i\neq j\neq k$ . Then $(n_i; v_i)\oplus (n_j; v_j) = (n_i; v_i)\oplus (n_k; v_k)$ if and only if $(n_j; v_j) = (n_k; v_k)$.
\el
\bp
Suppose $(n_j; v_j) = (n_k; v_k)$, then 
$$(n_i; v_i)\oplus (n_j; v_j) = (n_j\wedge n_i; v_i\vee v_j) = (n_k\wedge n_i; v_i\vee v_k) = (n_i; v_i)\oplus(n_k; v_k).$$
Conversely suppose that
$(n_i; v_i)\oplus (n_j; v_j) = (n_i; v_i)\oplus(n_k; v_k)$.
Then
$$(n_j\wedge n_i; v_i\vee v_j) = (n_i\wedge n_j; v_j\vee v_i) = (n_k\wedge n_i; v_i\vee v_k) = (n_i\wedge n_k; v_k\vee v_i).$$
 Since $v_j\leq n_i$ and $v_i\leq n_j,\, (n_j; v_j)(n_i; v_i)=(n_j; v_j)$ and 
\begin{eqnarray*}
(n_j; v_j)(v_k; n_k) &=& (n_j; v_j)(v_i; n_i)(v_k; n_k)\\
                     &=& (n_j; v_j)(v_i\vee v_k; n_i\wedge n_k)\\
										 &=& (n_j; v_j)(v_j\vee v_i; n_j\wedge n_i)\\
										 &=& (n_j; v_j)(v_j; n_j)(v_i; n_i)\\
										 &=& (1;0)(v_i; n_i)\\
										 &=& (1;0)										
\end{eqnarray*}
therefore, $S((n_j; v_j)(v_j; n_j)) = \left\{(1;0)\right\}$ and so $v_j\leq v_k$ and
\begin{eqnarray*}
(v_k; n_k)(n_j; v_j) &=& (v_k; n_k)(v_i; n_i)(n_j; v_j)\\
                     &=& (v_k\vee v_i; n_i\wedge n_k)(n_j; v_j)\\
										 &=& (v_i\vee v_j; n_j\wedge n_i)(n_j; v_j)\\
										 &=& (v_i; n_i)(v_j; n_j)(n_j; v_j)\\
										 &=& (v_i; n_i)(1;0)\\
										 &=& (1;0)										
\end{eqnarray*}
Therefore, $S((v_k; n_k),(n_j; v_j)) = \left\{(1;0)\right\}$ and so $n_j\leq n_k$.
Interchanging $(n_k; v_k)$ and $(n_j; v_j)$, $n_j\leq n_k$, thus $n_j = n_k$ and $v_j = v_k$. That is, $(n_j; v_j) = (n_k; v_k)$.
\ep

\bc
Let $(n_i; v_i), (n_j; v_j) \in E_{P(L)}$ with $v_j\leq n_i$ for $i\neq j$. Then $$(n_i; v_i)\oplus(n_j; v_j) = (0;1) \mbox{ if and only if }(n_j; v_j) = (v_i; n_i).$$
\ec

\bp
By Lemma.\ref{4a} we have
$$S((n_i; v_i),(v_i; n_i)) = S((v_i; n_i),(n_i; v_i)) = \left\{(1;0)\right\}.$$
 Therefore
$$S((n_i; v_i),(v_i; n_i))\cap S((v_i; n_i),(n_i; v_i)) = \left\{(1;0)\right\}.$$
ie., $((v_i\vee n_i; n_i\wedge v_i)) = (1;0)$, and thus 
$$(n_i; v_i)\oplus(v_i; n_i) = (n_i\wedge v_i; v_i\vee n_i) = (0;1).$$
Conversely, suppose $(n_i; v_i)\oplus(n_j; v_j) = (1;0)$, since
$(n_i; v_i)\oplus(v_i; n_i) = (1;0)$ it follows that
$(n_j; v_j) = (v_i; n_i)$. 
\ep

\bl
Let $E_{P(L)}$ be the biordered set with $v_i\leq n_j$, $v_j\leq n_i$, $v_i\leq n_k$, $v_k\leq n_i$, $v_j\leq n_k$, $v_k\leq n_j$ for $i\neq j\neq k$. Then for elements $(n_i; v_i), (n_j;v_j), (n_k; v_k)$, $i, j, k = 1,2,\hdots, N$ with $i\neq j\neq k$ in $E_{P(L)}$, the collection $\left\{v_1, v_2, \hdots, v_N\right\}$ are independent elements in the lattice $L$.
\el

\bp 
We have the set $\left\{(n_i; v_i)\colon i = 1,2,\hdots, N\right\}$ in $E_{P(L)}$ so that the elements $v_1, v_2, \hdots, v_N$ are in the complemented modular lattice $L$. 
Since $v_i\leq n_k$ and $v_j\leq n_k$ for $i\neq j$ in $E_{P(L)}$, we have $(v_i \vee v_j)\leq n_k$ for $i\neq j\neq k$, then
$$(v_i\vee v_j)\wedge v_k \leq n_k\wedge v_k = 0.$$
Hence the collection $\left\{v_i \colon i= 1,2 \hdots, N\right\}$ are independent.
\ep

Now recall from \cite{Akhi}, that given a regular ring $R$ whose biordered  set of idempotents $E_R$ then the principal biorder ideals $\omega^l(e)\,[\omega^r(e)], \,e\in E_R$ 
is the complemented modular lattice $\Omega_L\,[\Omega_R]$. Also we have the following theorem

\bt (\cite{Akhi}, Theorem 6) Let $R$ be a regular ring with idempotents $e_1,e_2,\cdots, e_n$ such that $M(e_1,e_j)=\{0\}$ for $i\neq j, \, d_l(e_1,e_j)\leq 3$ and $e_1+e_2+\cdots +e_n=1.$ 
Then a maximal system $\omega^l(e_1),\cdots \omega^i(e_n)$ in $\Omega_L$ forms a homogeneous basis of rank $n$ and the complemented modular lattice is of order $n$.
\et

It is easy to observe  that $E_{P(L)}$ has elements $\left\{(n_i;v_i)\colon i = 1,2, \hdots, N\right\}$ satisfying the following properties:

\begin{enumerate}
	\item $v_i \leq n_j$ for $i\neq j$
	\item $(n_1; v_1)\oplus (n_2; v_2)\oplus \hdots \oplus(n_N; v_N) = (0;1)$
	\item $d_l((n_i; v_i),(n_j; v_j)) = 3$ for $i \neq j$.
\end{enumerate}
and this subset of $E_{P(L)}$ is again a biordered set and we denote this biordered set as $E_{P(L)}^0$\index{$E_{P(L)}^0$}. Let $L^0$ be the complemented modular lattice 
with $E_{P(L^0)}=E_{P(L)}^0$ then $L^0$ is of order $N$.

The next lemma gives a biorder isomorphism between the biordered set of idempotents in the ring $R$ and the biordered set $E_{P(\Omega_L)}$.
	\bl
	Every idempotent $e$ in a ring $R$ is associated with a pair $(\omega^l(1-e); \omega^l(e))$ of complementary biorder ideals in $\Omega_L$. The map $\epsilon\colon E_R\rightarrow E_{P(\Omega_L)}$ defined by $\epsilon(e) = (\omega^l(1-e); \omega^l(e))$ is a biorder isomorphism.
	\el
	\bp
	For each $e \in E_R$, $(\omega^l(1-e); \omega^l(e))$ is a complementary pair in the lattice $\Omega_L$ and the mapping
	$$\epsilon\colon e\rightarrow (\omega^l(1-e); \omega^l(e))\,\,\mbox{for all } e \in E_R$$
	is a map of $E_R$ into $E_{P(\Omega_l)}$.
	The map $\epsilon$ is clearly injective. It follows from the definition of biordered set (\cite{Namb}, Definition 1) and the equation.1 in Theorem(1) \cite{pastjin} that the map 
	$\epsilon$ preserve basic products and hence $\epsilon\colon E_R\rightarrow E_{P(\Omega_L)}$ is a biorder isomorphism. Also, it can be easily seen that this map $\epsilon$ is a regular bimorphism.
	\ep
	Thus we have $E_{P(\Omega_L)}$ and $E_R$ are biorder isomorphic. Now we show that there exists elements $e_1, e_2, \hdots, e_N$ in $E_R$ satisfying all the conditions of 
	$E_{P(\Omega_L)}^0$. Consider $E_{P(\Omega_L)}^0$, then there are elements $((\oml(1-e_i); \oml(e_i)\colon i= 1,2, \hdots, N)$ such that
	\begin{enumerate}
		\item $\oml(1-e_i)\leq \oml(e_j)$ for $i \neq j$
		\item $(\oml(1-e_1); \oml(e_1))\oplus (\oml(1-e_2); \oml(e_2))\oplus \hdots \oplus (\oml(1-e_N); \oml(e_N)) = (0;1)$
		\item $d_l((\oml(1-e_i); \oml(e_i)), (\oml(1-e_j); \oml(e_j)))= 3$
	\end{enumerate} 
	Since $E_{P(\Omega_L)}$ and $E_R$ are biorder isomorphic, for each $((\oml(1-e_i); \oml(e_i)\colon i= 1,2, \hdots, N)$, there exists elements $e_1, e_2, \hdots, e_N$ such that
	\begin{enumerate}
		\item $\oml(1-e_i) \leq \oml(e_j)$ for $i\neq j$ implies $(\oml(1-e_i); \oml(e_i))(\oml(1-e_j); \oml(e_j)) = \left\{(0;1)\right\}$.
		But $\oml(1-e_i) \leq \oml(e_j)$ implies $1-e_i \ml\oml e_j$ thus $e_i \ml\omr (1-e_j)$ and so $e_ie_j = 0$ for $i \neq j$.\\
		
		\item Since $(\oml(1-e_1); \oml(e_1)) \oplus \hdots \oplus (\oml(1-e_N); \oml(e_N)) = (0; \omega^l(1))$
		which implies $\oml(e_1)\vee \oml(e_2) \vee \hdots \vee \oml(e_N) = \oml(1)$. But we have from  (\cite{Akhi}, Lemma 9) $e_ie_j = 0$ for $i \neq j$ implies $\oml(e_1)\vee \oml(e_2)\vee \hdots \vee\oml(e_N) = \oml(e_1+e_2+ \hdots +e_N) = \oml(1)$ and hence $e_1+e_2+ \hdots +e_N = 1$.
		\item $d_l((\oml(1-e_i); \oml(e_i)), (\oml(1-e_j); \oml(e_j)))= 3$ implies there exists elements $(\oml(1-e_i); \oml(e_i))\mathcal{L}$ $(\oml(1-e_h); \oml(e_h)) \mathcal{R}$ $(\oml(1-e_k); \oml(e_k))\mathcal{L}(\oml(1-e_j); \oml(e_j))$. Therefore, by definition of $\mathcal{L}$ and $\mathcal{R}$, $\oml(e_i) = \oml(e_h)$ and $\oml(1-e_h) = \oml(1-e_k)$ and 
		$\oml(e_k) = \oml(e_j)$. But $\oml(1-e_h) = \oml(1-e_k)$ implies $\omr(e_h) = \omr(e_k)$. Thus we get $e_i \mathcal{L}e_h\mathcal{R}e_k \mathcal{L} e_j$ and hence $d_l(e_i, e_j) = 3$.
		\end{enumerate}
		Since $E_{P(\Omega_L)}^0$ is a biorder subset of $E_{P(\Omega_L)}$ and corresponding to each element in $E_{P(\Omega_L)}$ there exists elements in $E_R$ satisfying all the conditions of $E_{P(\Omega_L)}^0$ as shown above, we have $E_{P(\Omega_L)}^0$ and $E_{R}^0$ (where $E_{R}^0$ is the biordered set similar to $E_{P(L)}^0$) are biorder isomorphic.
	
	Thus it can be seen that for given a complemented modular lattice $L$ there exists a complemnted modular lattice $L^0$ admiting a biordered subset $E_{P(L)}^0$ of the biordered set 
	$E_{P(L)}$ consisting of $N$ elements and hence $L^0$ admits a homogeneous basis of order $N$. 
	Thus analogous to von-Neumann's coordinatization theorem, we have the following theorem:
	\bt
	Let $L$ be a complemented modular lattice admitting a biordered subset with at least $4$ elements, having the following properties:
	\begin{enumerate}
		\item $v_i \leq n_j$ for $i\neq j$
	\item $(n_1; v_1)\oplus (n_2; v_2)\oplus \hdots \oplus(n_N; v_N) = (0;1)$
	\item $d_l((n_i; v_i),(n_j; v_j)) = 3$ for $i \neq j$,
	\end{enumerate} 
	 then there exists a von Neumann regular ring $R$ such that $L$ is isomorphic to the lattice of all principal left ideals of $R$.
	\et

\end{document}